\documentclass[10pt]{amsart}

\usepackage{vmargin}
\usepackage{amsmath}
\usepackage{amssymb}
\usepackage{amsfonts}
\usepackage{enumerate}
\usepackage{stmaryrd,color}

\numberwithin{equation}{section}

\newcommand{\bC}{{\mathbb{C}}}

\newcommand{\bN}{{\mathbb{N}}}

\newcommand{\bR}{{\mathbb{R}}}

\newcommand{\bZ}{{\mathbb{Z}}}

  \newcommand{\B}{{\mathcal{B}}}

  \newcommand{\E}{{\mathcal{E}}}

\renewcommand{\H}{{\mathcal{H}}}

  \newcommand{\M}{{\mathcal{M}}}
  \newcommand{\N}{{\mathcal{N}}}

\renewcommand{\phi}{\varphi}
\newcommand{\upchi}{{\raise.35ex\hbox{\ensuremath{\chi}}}}
\newcommand{\eps}{\varepsilon}

\renewcommand{\leq}{\leqslant}
\renewcommand{\geq}{\geqslant}

\newcommand{\norm}[1]{\left\| #1 \right\|} \newcommand{\md}[1]{\left|
  #1 \right|} \newcommand{\p}[1]{\left( #1 \right)}

\newtheorem{thm}{Theorem}[section] \newtheorem{defi}[thm]{Definition}
\newtheorem{prop}[thm]{Proposition} \newtheorem{cor}[thm]{Corollary}
\newtheorem{lemma}[thm]{Lemma} \newtheorem{rk}[thm]{Remark}

\begin{document}
\title{Revisiting the Marcinkiewicz theorem for non commutative maximal functions} \date{}
\author{L\'eonard Cadilhac}\address{Institut de Math\'ematiques de Jussieu, Sorbonne Universit\'e - Paris, France}\email{cadilhac@imj-prg.fr} \author{\'Eric
  Ricard} \address{UNICAEN, CNRS, LMNO, 14000 Caen,
  France} \email{eric.ricard@unicaen.fr}

\maketitle

\begin{abstract}
We give an alternative proof of a Marcinkiewicz interpolation theorem for non commutative maximal functions and positive maps and refine earlier versions of the statement. The main novelty is that it provides
a substitute for the maximal function of a martingale in $L_p$,
$1<p\leq \infty$, losing very little on numerical constants. For non positive maps, the above mentioned theorem fails but we can still obtain some interpolation results by weakening the maximal norm that we consider.
\\

\smallskip
\noindent \textbf{Keywords.} non commutative integration, real interpolation, maximal inequalities, weak type inequalities.
\\

\smallskip
\noindent \textbf{2020 Mathematics Subject Classification.} 46L51, 46B70, 46L52.

\end{abstract}
\section{Introduction}

In classical analysis, a maximal function $M$ is simply the supremum of a given family of functions $(f_i)_{i\in I}$,
$$
M := \sup_{i\in I} \md{f_i}.
$$
Obtaining bounds for certain meaningful maximal functions is a fundamental theme in harmonic analysis, with a variety of applications ranging from Lebesgue's differentiation theorem, to ergodic theory and the convergence of Fourier series. When trying to adapt this notion to non commutative analysis, where functions are replaced by operators, an immediate difficulty appears: in general, a family of operators does not admit a supremum. Nonetheless, alternative ways to formulate maximal inequalities exist and can be traced back at least to the 70's and the early days of non commutative martingale and ergodic theory (\cite{Cuc71}, \cite{Lan76}, \cite{Yea77}). A systematic approach to non commutative maximal functions was initiated two decades ago in works of Pisier \cite{P} and Junge \cite{J2} and has been widely adopted since.

In this paper, we revisit a result of Junge and Xu \cite{JX2} that extends Marcinkiewicz interpolation in the following way: let $p_0 < p <p_1$, if a non commutative positive maximal operator is weakly bounded on $L_{p_0}$ and $L_{p_1}$, then it is strongly bounded on $L_p$. Interpolation generally transfers well from classical to non commutative function spaces so this may appear to be a routine result but - details will be given shortly - weak and strong type have here been defined in different ways and independently (weak type in \cite{Cuc71} and \cite{Yea77} and strong type in \cite{P}) making this theorem, if not surprising, quite remarkable. Our first goal is to present a new, very simple proof of this theorem and then to explore its potential extensions to operators that are not necessarily positive. A motivation discussed along the way is to compare different definitions of maximal norms and how real interpolation might apply to them. Note that simplifications and refinements of Junge and Xu's original result can already be found in \cite{BCO}, \cite{DirkJOT}. In particular, Junge and Xu always assume strong type $(p_1,p_1)$ and the result at the beginning of the paragraph was only proved later by Dirksen \cite{DirkJOT}.

Let us now recall the definition of maximal norms as introduced by Pisier and Junge. Let $(\N,\tau)$ 
be a von Neumann algebra equipped with a normal semifinite
faithful trace. 
For $0<p\leq\infty$, the space $L_p(\N;\ell_\infty)$ consists of all
sequences of operators $x=(x_n)_{n\geq0}$ for which there exist operators 
$a,b\in L_{2p}(\N)$ and contractions $(u_n)_{n\geq0}$ in $\N$ such that
\begin{equation}\label{eq:maxnorm}
    \forall n\geq1,\qquad x_n =au_nb.
    \end{equation}
    We set $$\|x\|_{L_p(\N;\ell_\infty)}=\inf \|a\|_{2p} \|b\|_{2p},$$
where the infimum runs over all decompositions as above. Then $L_p(\N;\ell_\infty)$ is a Banach space for $p\geq 1$ and a quasi-Banach for $p<1$. We refer 
to \cite{J2} for further details.
When $x$ consists of positive  operators, we have a simpler description
 $$\|x\|_{L_p(\N;\ell_\infty)}=\inf \{\|a\|_{p} \;| \forall n\geq 0\;,
0\leq x_n\leq a,\}.$$ 
This explains why $\|x\|_{L_p(\N;\ell_\infty)}$
plays the role of the norm of the maximal function. It is proved in \cite{JX} that this infimum is achieved but it may depend on $p$ if $x$ belongs to $L_p(\N;\ell_\infty)$ for multiple values of $p$. Advantages of these definitions include: a good behaviour with respect to complex interpolation \cite{J2, XW}, a nice duality theory and the ability to carry out in the non commutative setting several classical applications of maximal functions. A drawback is that these spaces do not form a real interpolation scale \cite{JX3}.

We will also need a weak-version that is commonly used in non commutative
analysis. The space $\Lambda_{p,\infty}(\N;\ell_\infty)$ consists of all
sequences of operators $x=(x_n)_{n\geq0}$ for which the following quantity is finite
$$\|x\|_{\Lambda_{p,\infty}(\N;\ell_\infty)}= \sup_{\lambda>0} \lambda \inf_{e} 
\{\big(\tau(1-e)\big)^{1/p} \; |\; \forall n\geq 0,\,\|ex_ne\| \leq \lambda\},$$
where the infimum runs over all self-adjoint projections in $\N$. 
Thus, if
$x\in  \Lambda_{p,\infty}(\N;\ell_\infty)$ with $\|x\|_{\Lambda_{p,\infty}(\N;\ell_\infty)}<C$ if
 for all $\lambda>0$, there is a projection 
$e\in\N$ with
\begin{equation}\label{normefaible}
  \tau(1-e)\leq \frac {C^p}{\lambda^p},\, \qquad 
\forall n\geq 0,\,\|ex_ne\| \leq \lambda.\end{equation}
$\Lambda_{p,\infty}(\N;\ell_\infty)$ is a quasi-Banach space.

At this stage, we point out that it does not match with a more natural 
definition as for commutative measure spaces commonly denoted by $L_{p,\infty}(\N;\ell_\infty)$. If $x$ consists of positive operators, $x\in \Lambda_{p,\infty}(\N;\ell_\infty)$ does not implies that there exists $a\in L_{p,\infty}(\N)$ such that
$0\leq x_n\leq a$. This requirement would provide better constants for interpolation (see Remark \ref{strongweak}) but is too strong in practice to formulate Doob's weak type inequality for non commutative martingales or a weak type inequality for ergodic theory. In the following, $\M$ denotes another von Neumann algebra equipped with a n.s.f. trace. 

\begin{defi}
Let $S=(S_n)_{n\geq 0}$ where $S_n:L_p(\M)\to L_0(\N)$ be a sequence of maps
  \begin{itemize}
  \item $S$ is said to be of strong type $(p,p)$ with constant $C$ if
    $S$ is bounded from $L_p(\M)$ to $L_p(\N;\ell_\infty)$ with
    constant $C$.
 \item $S$ is said to be
of weak type $(p,p)$ if $S$ is bounded from $L_p(\M)$ to $\Lambda_{p,\infty}(\N;\ell_\infty)$.\\
 More precisely, we say that $S$ is of weak type $(p,p)$ with constant $C$, if 
for any $x\in L_p(\M)$ and $\lambda>0$, there 
exists a projection $e\in \N$ such that 
\end{itemize}
\begin{equation}\label{weaktype}
\tau(1-e)\leq \frac {C^p\|x\|_p^p}{\lambda^p},\, \qquad 
\forall n\geq 0,\,\|eS_n(x)e\| \leq \lambda.
\end{equation}
\end{defi}

Note that when $p=\infty$, the weak and strong types are equivalent and simply mean that the family $(S_n)$ is uniformly bounded from $\M$ to $\N$.

The  family $(S_n)$ is of weak type $(p,p)$ with constant $C$ iff $(S_n):L_p(\M)\to \Lambda_{p,\infty}(\N;\ell_\infty)$ with norm $C$.

\smallskip

We are now able to state the Marcinkiewicz interpolation theorem 
for non commutative maximal function. 

\begin{thm}[Junge, Xu, Dirksen]\label{thm:JXD}
 Let $1 \leq p_0 < p_1 \leq \infty$. Let $S = (S_n)_{n\geq 0}$be a sequence of positive linear maps from $L_{p_0}(\M)+ L_{p_1}(\M)$ to $L_0(\N)$. Assume that $S$ is of weak type $(p_0,p_0)$ with constant $C_0$ and of weak type $(p_1,p_1)$ with constant $C_1$. Let $\theta \in (0,1)$ and $p$ be determined by $1/p = (1-\theta)/p_0 + \theta/p_1$. Then $S$ is of strong type $(p,p)$ with constant less than
 $$
 CC_0^{1-\theta}C_1^{\theta} \alpha_{\theta,p_0,p_1}^2,
 $$
 where $C$ is a universal constant and 
 $$
  \alpha_{\theta,p_0,p_1} = \Big(\frac 1 {p}-\frac 1 p_0\Big)^{-1} + \Big(\frac 1 {p_1}-\frac 1 p\Big)^{-1}.
  $$
\end{thm}

Compared to earlier approaches (\cite{JX2},\cite{BCO},\cite{DirkTAMS},\cite{DirkJOT}) ours has mainly two advantages: first, it identifies clearly the non commutative part of the proof which is reduced to a single lemma (Lemma \ref{majdiag}), the rest of the argument consists in manipulations of singular values ; second, it allows to construct an element in $L_p(\N)$, playing the role of a maximal function. Keeping the notations of Theorem \ref{thm:JXD}, we have

\begin{thm}
Let $x \in L_{p_0}(\M) \cap L_{p_1}(\M)$. Then, there exists $a \in (L_1 + L_{\infty})(\N)^+$ and contractions $(u_n)_{n\geq 0}$ such that 
$$
\forall n\geq 0,\quad S_n(x) = au_na.
$$
Moreover, for any $p \in (p_0,p_1)$
$$
\norm{a}_{2p}^2 \leq C_0^{1-\theta}C_1^{\theta} \alpha_{\theta, p_0,p_1}^{2} \big(\ln  \alpha_{\theta, p_0,p_1}\big)^2,
$$
where $p$ and $\theta$ still verify $1/p = (1-\theta)/p_0 + \theta/p_1$.
\end{thm}

As previously mentioned, in general, the minimizing factorization $(x_n) = (ay_nb)$ in \eqref{eq:maxnorm} may depend on $p$. This is an important difference with the commutative case and can be considered to be at the root of the failure of real interpolation for non commutative maximal norms. We show that in this particular case, if we are willing to lose an exponent in the constants, this phenomenon cannot occur and a concrete maximal function can be exhibited. In particular, this applies to $(S_n)_{n\geq 0}$ a family of conditional expectation to define a maximal function for non commutative martingales. We state the results for sequences of maps for convenience, they hold as well for general families
$S=(S_i)_{i\in I}$.

In the rest of the paper, we explore variants and limitations of this theorem. First, in the remainder of Section \ref{sec:int}, we show that the method of proof extends to some asymmetric versions of the maximal norm, already considered in \cite{J2}, \cite{GJP16} and \cite{XW}. For those, the range of values of $p$ for which the theorem applies has to be restricted (see Theorem \ref{thm:asym} and Proposition \ref{prop:noasym}). 

We begin Section \ref{sec:beyondpos} by noting that without the positivity hypothesis on the maps $(S_n)$, the theorem fails even with $\M = \mathbb{C}$. We introduce weaker maximal quasi-norms taking their inspiration in the definition of weak type rather than strong type. We show that they form a real interpolation scale and do give rise to Banach spaces for some parameters. They seem to be the most well-adapted kind of maximal inequalities beyond the positive case. 
They can be useful if one's goal is only to study questions of pointwise convergence in tracial von Neumann algebras.
We also show that one of them exactly corresponds to the real interpolated spaces for the couple $\big(L_p(\M;\ell_\infty^c), L_\infty(\M;\ell_\infty^c)\big)$ of one-sided maximal function spaces. We conclude with some basic bounds related to asymmetric factorizations of the form: 
\begin{equation}
    x_n = ay_n + z_nb \quad \text{with} \quad (y_n),(z_n) \in L_{\infty}(\M;\ell_\infty)\ \ \text{and}\ a,b\in L_p(\M).
\end{equation}  This is close in spirit to what Junge and Xu did in \cite{JX2}
for symmetric factorizations. Unfortunately, one can not improve them  using  of the machinery of interpolation as the real method fails for $(L_p(\M;\ell_\infty^c))_{p\geq 1}$.

\section{Preliminaries}

We assume that the reader is familiar with
non commutative integration. This section briefly recalls some standard notations and definitions (see also \cite{P}), and presents a simple decomposition lemma. 

\subsection{Non commutative integration}

We will denote by $(\M,\tau_\M)$ or $(\N,\tau_\N)$ non commutative measure spaces, meaning that $\M$ and $\N$ are von Neumann algebras and $\tau_\M$ (resp. $\tau_\N$) is a semifinite normal faithful trace on $\M$ (resp. on $\N$). In practice there will be no ambiguity on the trace used and we will write $\tau$ instead of $\tau_\N$ or $\tau_\M$. The space of $\tau$-measurable operators affiliated with $\M$ is denoted by $L_0(\M)$ and the non commutative $L_p$-spaces associated with $(\M,\tau)$ by $L_p(\M)$, $p\in (0,\infty]$. For $x \in L_0(\M)$, $\mu(x)$ designates the singular value function of $x$. 

The first step of the main proofs of this paper is to obtain estimates for finite projections. The following simple lemma is essential to extend those estimates to more general operators. 

\begin{lemma}\label{decdya}
Let $p\in (0,\infty)$ and $x\in L_p^+(\M)$. There are finite projections $(r_n)_{n\in \bZ} \in \M$ so 
that $x=\sum_{n\in \bZ} 2^{-n} r_n \in L_p(\M)$ and for all $\alpha>0$ 
$$
\sum_{n\in \bZ} 2^{-n\alpha} 1_{[0,\tau(r_n)]}\leq \frac 1
{1-2^{-\alpha}}\mu(x^\alpha) \quad \textrm{and} \quad \sum_{n\in\bZ} (\md{n}+1) 2^{-n} 1_{[0,\tau(r_n)]} \leq C \mu(x(\textrm{ln}(\md{x})+1))
$$
\end{lemma}
\begin{proof}
Actually, this is a commutative result.  Let $\phi_n$ be the indicator function on $\bR$ of the set of reals 
whose $-n^{\mathrm{th}}$-digit in base 2 is 1. Clearly 
for all $t\geq 0$, $t=\sum_{n\in \bZ} 2^{-n} \phi_n(t)$.
Set $r_n=\phi_n(x)$, it is a finite projection as $x\geq 2^{-n} r_n$, 
thus $\tau(r_n)\leq 2^{np}\|x\|_p^p$. By Lebesgue's dominated convergence theorem $x=\sum_{n\in \bZ} 2^{-n} r_n$ holds in $L_p(\M)$.

Let $s>0$. If for all $k\in \bZ$, $\tau(r_k)<s$, then $\sum_{k\in \bZ}
2^{-k\alpha} 1_{[0,\tau(r_k)]}(s)=0$ and there is nothing to prove.
Otherwise set $k_s=\min \{ k \;|\; \tau(r_k)\geq s\}$, it is well
defined as $\tau(r_k)\leq 2^{kp}\|x\|_p^p$. Then $\sum_{k\in \bZ}
2^{-k\alpha} 1_{[0,\tau(r_k)]}(s)\leq \sum_{k\geq k_s}2^{-k\alpha}=
\frac {2^{-\alpha k_s}}{1-2^{-\alpha}}$. But $x^\alpha\geq
2^{-k_s\alpha} r_{k_s}$ (recall that they commute), thus
$\mu(x^\alpha)(s)\geq 2^{-k_s\alpha}\mu(r_{k_s})(s)=2^{-k_s\alpha}$. For the second inequality, follow the same proof and note that $\sum_{k\geq k_s} (\md{k} +1)2^{-k} \leq (\md{k_s}+1)2^{-k_s} \sum_{k\geq 0} (k+1)2^{-k}$.
\end{proof}

\begin{rk}
Lemma \ref{decdya} still holds for any $x\in L_0(\M)^+$ but in general the sum $\sum_{n\in\bZ} 2^{-n} r_n$ only converges in $L_0(\M)$.
\end{rk}

\subsection{Hardy-Littlewood majorization}
Let $f$ and $g$ be two non-increasing functions from $(0,\infty)$ to $\mathbb{R}^+$. We say that {\it $g$ majorizes $f$} and write $f \preceq g$ if
$$
\forall t>0,\qquad \int_0^t f \leq \int_0^t g.
$$
We will use the following properties. First, let $p\in [1,\infty]$ and $x,y \in L_p(\M)$, then 
$$
\mu(x) \preceq \mu(y) \Rightarrow \norm{x}_p \leq \norm{y}_p. 
$$
Second, assume that $x = \sum_{n\in \mathbb{Z}} x_n$, where the sum converges in $L_p(\M)$. Then 
\begin{equation}\label{eq:majsum}
\mu(x) \preceq \sum_{n\in\mathbb{Z}} \mu(x_n). 
\end{equation}
One way to justify this inequality is to note that for any $t>0$, $x \mapsto \int_0^t \mu(x) = \norm{x}_{L_1 + tL_\infty}$ is a norm and use the triangle inequality. 

\section{Interpolation for positive maps}\label{sec:int}

 Let $(\M,\tau)$ be a semifinite von Neumann algebra equipped with a n.s.f. trace. 

\subsection{Majorization and factorization}

We start by two easy lemmas that capture the non commutative aspects
of the proof of our main theorem.

\begin{lemma}\label{majdiag}
 Let $N\in \bN$ and $(q_k)_{|k|\leq N}$ be a sequence of disjoint projections in $\M$
with  $e=\sum_{|k|\leq N} q_k$. Let  $(d_k)_{|k|\leq N}$ be a 
 sequence of strictly positive reals. Then for any $x \in L_0(\M,\tau)^+$:
$$0\leq exe\leq \Big(\sum_{|k|\leq N} 1/d_k\Big)\sum_{|k|\leq N} d_k q_k x q_k.$$
\end{lemma}
\begin{proof}
 The matrix $(1/\sqrt{ d_id_j})_{-N\leq i,j\leq N}$ corresponds to a
 rank one positive operator with norm $C=\sum_{|k|\leq N} 1/d_k$. Thus
 the matrix $(C\delta_{i,j}-1/\sqrt{ d_id_j})_{-N\leq i,j\leq N}$ is
 positive, hence conjugating by $(\sqrt{d_i}\delta_{i,j})$,
 $(Cd_i\delta_{i,j} -1)_{-N\leq i,j\leq N}$ also is. We can find
 $T \in \bN$ and families of complex numbers $(a_{i,t})_{\md{i}\leq N, t\leq T}$ such that 
$Cd_i\delta_{i,j} -1=\sum_{t=1}^T  \overline{a_{i,t}}a_{j,t}$. Hence 
$$C\sum_{|k|\leq N} d_k q_k x q_k-exe= \sum_{t=1}^T  \Big(\sum_{|i|\leq N}{a_{i,t}}q_i\Big)^* x \Big(\sum_{|j|\leq N}{a_{j,t}}q_i\Big)\geq 0.$$
\end{proof}

The following is standard by using polar decompositions in $\M\overline{\otimes} B(\ell_2(\bZ))$:
\begin{lemma}\label{fact}
  Let $(a_i)_{i\in\bZ}$ and $(b_i)_{i\in\bZ}$ be sequences of elements in 
$L_r$ and $L_s$ such that 
$\sum_{i\in \bZ} a_ia_i^*\in L_{r/2}$ and $\sum_{i\in \bZ} b_i^*b_i\in L_{s/2}$ with $r,s>0$. 
For any sequence of contractions $(u_i)_{i\in\bZ}$ in $\M$, there is a contraction $u\in \M$ such that, in $L_{rs/(r+s)}$:
$$\sum_{i\in\bZ} a_i u_i b_i =  \Big(\sum_{i\in \bZ} a_ia_i^*\Big)^{1/2} u \Big(\sum_{i\in \bZ} b_i^*b_i\Big)^{1/2}.$$
\end{lemma}

\subsection{The Marcinkiewicz theorem}

The first version of the Marcinkiewicz interpolation for maximal 
functions is given by \cite[Theorem 3.1]{JX2}: 
assuming weak type $(p_0,p_0)$ and strong type $(p_1,p_1)$, 
the authors obtain strong type $(p,p)$ for $p_1 < p < p_0$. This has been extended in several directions in \cite{DirkTAMS, BCO}. 
A very satisfactory statement was obtained by Dirksen in \cite{DirkJOT}. Our approach is similar to Dirksen's, the novelties are in the second statement and in the simplicity of 
the proof which is ``almost commutative'' apart from Lemmas \ref{majdiag} and \ref{fact}. Actually, a careful look at \cite{DirkJOT}, fixing
$p_0<r<s<p_1$, one can find $a$ depending only on $r,s$ such that 
\eqref{mainest} hold for $r<p<s$.
 
\begin{thm}\label{Mar}
  Let $1\leq p_0< p_1\leq \infty$ and $S=(S_n)$ be a sequence of
  linear positive maps from $L_{p_0}(\M)+ L_{p_1}(\M)$ to
  $L_0(\N)$. Assume that $S$ if of weak type $(p_0,p_0)$ and
  $(p_1,p_1)$ with constants $C_0$ and $C_1$. Let $p \in (p_0,p_1)$ and $\theta \in (0,1)$ such that $\frac 1 {p}= \frac{1-\theta}{p_0}+\frac \theta{p_1}$. Set 
  $$
  \alpha_{\theta,p_0,p_1} = \Big(\frac 1 {p}-\frac 1 p_0\Big)^{-1} + \Big(\frac 1 {p_1}-\frac 1 p\Big)^{-1}.
  $$
Then $S$ is of strong type $(p,p)$ with constant 
$$C_{\theta,p_0,p_1} \leq C C_0^{1-\theta} C_1^\theta  \alpha_{\theta,p_0,p_1}^2.$$  
Moreover, there is a constant $C_{\theta,p_0,p_1}'$ such that for any $x\in
L_{p_0}(\M)\cap L_{p_1}(\M)$, there exist $a\in(L_{1}+L_\infty)^+$ , and contractions $u_n\in \N$ independent of $p$ or $\theta$ with 
\begin{eqnarray}\label{mainest}S_n(x)= a u_n a,\qquad  \|a\|_{2p}^2\leq C_{\theta,p_0,p_1}' \|x\|_{p}\end{eqnarray}
and.
$$
C_{\theta,p_0,p_1}' \leq CC_0^{1-\theta} C_1^\theta  \Big(\alpha_{\theta,p_0,p_1}\mathrm{ln}( \alpha_{\theta,p_0,p_1})\Big)^2.$$  
\end{thm}

To prove Theorem \ref{Mar}, we proceed in two steps. First, in Lemma \ref{basic}, assuming that $x$ is a projection, we give a construction of a majoring element based on Lemma \ref{majdiag}. Second, for a general $x$, we decompose $x$ into a dyadic combination of projections $ x = \sum 2^n x_n$, apply the first step to each projection and put them back together.

\smallskip

Scaling the trace, we can and will assume $C_0=C_1$ and by homogeneity $C_0=1$.

Assume $p_0 < p < p_1$ fixed. Let us introduce notations that will play an important role in the proof of Theorem \ref{Mar}. Set $I = \bZ$ if $p_1 \neq \infty$ and $I = \bZ_{\geq 0}$ if $p_1 = \infty$. Let $(d_k)_{k\in I}$ be a sequence of positive real numbers such that:
$$
C_d = \sum_{k\in I}\dfrac1{d_k} < \infty.
$$
Set $\tilde d_k=4C_d d_k2^{-k/p_0}$ if $k\geq 1$ and $\tilde d_k=4C_d d_k2^{-k/p_1}$ otherwise. Assume that 
\begin{equation}\label{eq:dktilde}
\sum_{k\in I} \tilde d_k 2^{k/p} < \infty. 
\end{equation}
For $s>0$, define the dilation operator $D_s$ by 
\begin{align*}
D_s\colon L_0(\bR^+) & \longrightarrow L_0(\bR^+)\\
f &\longmapsto [t \mapsto f(t/s)].
\end{align*}
The connection between dilation operators and Marcinkiewicz interpolation was made explicit by Boyd in \cite{Boy69}. A closely related way to approach interpolation of weak type inequalities was earlier developed by Calder\'on \cite{Cal66}. We found Boyd's formulation to be more convenient in this paper but both are essentially equivalent (as shown for example in \cite{Mon96}). 

\begin{lemma}\label{basic}
  Let $r\in \M$ be a finite projection. There exists an element  
   $z\in L_p(\N)$ (depending only on $r$ and the choice of $(d_k)_{k\in \bZ}$) such that
   $$\forall \;n\geq 0, \quad 0\leq S_n(r)\leq z, \qquad \mu(z)\preceq \sum_{k\in I} \tilde d_kD_{2^k}(\mu(r)).$$
\end{lemma}

\begin{proof}
Let  $r\in \M$ with $\tau(r)=t$. 

The map $S$ is of weak type $(p_0, p_0)$ with constant $1$. Hence, using \eqref{weaktype} for
$\lambda=2^{-(k-2)/p_0}$ and $k\geq 1$, we obtain projections
$(e_k)_{k\geq 1}$ such that
$$ \forall k\geq 1\; \forall \;n\geq0, \quad \| e_k S_n(r)e_k\| \leq   2^{-(k-2)/p_0}, \quad
\tau(1-e_k)\leq 2^{k-2} t.$$
Similarly the weak type $(p_1,p_1)$ of $S$ gives  projections $(e_k)_{k\leq 0}$ such that
$$ \forall \;n\geq0, \quad \| e_k S_n(r)e_k\| \leq    2^{-(k-2)/p_1}, \quad
\tau(1-e_k)\leq 2^{k-2} t.$$

Considering $1- \vee_{i\leq k}(1-e_i)$ instead of $e_i$, we get a decreasing family of
projections such that
$$ \forall \;n\geq0, \quad \| e_k S_n(r)e_k\| \leq  2^{2/p_0}   2^{-k/p_0} 1_{k>0}+ 2^{2/p_1}  2^{-k/p_1} 1_{k\leq 0}, \quad
\tau(1-e_k)\leq 2^{k-1}t.$$

Set $q_k=e_{k}-e_{k+1}$, $E_N=\sum_{|k|\leq N} q_k$, we have $\sum_{k\in \bZ} q_k=1$ and $\tau(q_k)\leq 2^{k}t$. Note that if $p_1=\infty$, we can take $e_k=1$ for $k\leq 0$
so that $q_k=0$ for $k<0$.

By Lemma \ref{majdiag}, for all $n, N\geq0$
\begin{eqnarray*}  0\leq E_NS_n(r)E_N &\leq& C_d \sum_{|k|\leq N} d_k q_k S_n(r) q_k
\\ &\leq &4C_d \Big(  \sum_{k=1}^N d_k2^{-k/p_0} q_k +  
\sum_{k=0}^{N} d_{-k}2^{k/p_1} q_{-k}\Big)=z_N.
\end{eqnarray*}
By \eqref{eq:dktilde}, the increasing sequence $(z_N)_{N\geq 1}$ converges in
$L_{p}(\N)$ (using that $\tau(q_k)\leq 2^{k}t$). Denote by $z \in L_p(\N)$ its limit. Then
 $$\forall \;n\geq 0, \quad 0\leq S_n(r)\leq z=\sum_{k\in I}\tilde d_k q_k.$$ The inequality
$\tau(q_k)\leq 2^{k}t$ yields $\mu(q_k)\leq D_{2^k}(1_{[0,t]})$. Hence by \eqref{eq:majsum}
$$\mu(z)\preceq \sum_{k\in I} \tilde d_kD_{2^k}(\mu(r)),$$
with appropriate changes if $p_1=\infty$.
\end{proof}

\begin{proof}[Proof of Theorem \ref{Mar}]

  Writing any element $x\in L_{p}(\M)$ as a linear combination of 4 positive
  elements with coefficients $i^k$, each having norm less than $\|x\|_{p}$ and 
using Lemma \ref{fact}, it suffices to consider only $x\geq0$. This will only change the constant $C$.

Let $x\in L_{p}^+$, we use Lemma \ref{decdya} to write
$x=\sum_{m\in \bZ} 2^{-m} r_m$. We apply Lemma \ref{basic} to each
$r_m$ to obtain some element $z_m \in L_p(\M)$ with 
$$
\forall n\geq 0 \quad S_n(r_m) \leq z_m \quad \text{and} \quad 
\mu(z_m) \preceq \sum_{k\in I} \tilde d_k D_{2^k}(\mu(r_m)),
$$
where $I=\bZ$ if $p_1<\infty$ and $I=\bZ_{\geq 0}$ if $p_1=\infty$. For any finite subset $J\subset \bZ$:
\begin{equation}
\label{linear}\forall \;n\geq 0, \quad 0\leq S_n\p{\sum_{m\in J} 2^{-m}r_m} \leq z_J := \sum_{m\in J} 2^{-m} z_m.
\end{equation}
By \eqref{eq:majsum},
\begin{equation}\label{eq:majzI}
\mu\p{z_J} \preceq \sum_{m\in J} 2^{-m}\mu(z_m)
\preceq \sum_{m\in J}\sum_{k\in I} 2^{-m} \tilde d_kD_{2^k}(\mu(r_m))
\preceq \sum_{k\in I} \tilde d_kD_{2^k}\p{\sum_{m\in J} 2^{-m}\mu(r_m)}.
\end{equation}
By the construction of Lemma \ref{decdya}, $\sum_{m\in \bZ}2^{-m}\mu(r_m)\leq 2 \mu(x)$, thus
the sum $\sum 2^{-m}\mu(r_m)$ converges in $L_p(0,\infty)$. This combined with \eqref{eq:majzI} and \eqref{eq:dktilde} implies that the sum $\sum 2^{-m}z_m$ converges in $L_p(\N)$. Set $z = \sum_{m\in \bZ} 2^{-m}z_m$, we have
$$
\mu(z)\preceq \sum_{k\in I} \tilde d_kD_{2^k}\p{\sum_{m\in \bZ} 2^{-m}\mu(r_m)}
\leq 2 \sum_{k\in I}  \tilde d_kD_{2^k}(\mu(x)).
$$
In particular, 
$$
\| z\|_{p} \leq 2\sum_{k\in I}  \tilde d_k \|D_{2^k}(\mu(x))\|_p=
 2\|x\|_{p}\sum_{k\in I}  \tilde d_k 2^{k/p}.$$
Moreover, remark that by real interpolation, each $S_n$ is continuous on $L_p$ since it is of weak type $(p_0,p_0)$ and $(p_1,p_1)$, so 
$$
\forall n\geq 0 \quad S_n(x) = S_n\p{\sum_{m\in \bZ} 2^{-m}r_m} \leq \sum_{m\in\bZ} z_m = z.
$$

To get the first statement, we choose
 $d_k= 2^{k(1/p_0-1/p)/2}$ for $k\geq 1$ and
 $d_k= 2^{k(1/p_1-1/p)/2}$ for $k\leq 0$. We get that $\tilde d_k 2^{k/p}=4C_d/d_k$ for $k\in I$ where $C_d=\sum_{k\in I} 1/d_k$. Hence, computations give
$$\| z\|_{p}\leq   8\Big(\frac 1 {1-2^{(1/p-1/p_0)/2}}+\frac 2 
{1-2^{(1/p_1-1/p)/2}}\Big)^2\|x\|_{p};  $$ 
when $p_1=\infty$, we get only 1 for the second term in the sum.

To prove the second statement we choose $d_k={|k|(\ln|k|)^2+1}$, thus
$z$ is independent of $p$ and 
\begin{equation}\label{majz}\| z\|_{p}\leq   C \|x\|_{p}\sum_{k\geq 0}  (k(\ln|k|)^2+1) (2^{-k(1/p_0-1/p)}+
2^{-k(1/p-1/p_1)}).\end{equation}
 
 Up to an absolute factor, $\sum_{k\geq 0}  (k(\ln|k|)^2+1) 2^{-kx}$, $x>0$ is controlled by 
 $$\int_2^\infty 2^{-tx} t |\ln t|^2dt \lesssim  \frac {|\ln x|^2+1} {x^2}.$$ 
This yields the estimate.


\end{proof}

\begin{rk}\label{casinfini}{\rm When $p_1=\infty$, the constant $C_{\theta, p_0,\infty}$ actually remains bounded when $\theta\to 1$.
 One can also see from the proof that $z$ and $a$ are
 bounded when $x$ is bounded. Indeed, by construction the family $(z_m)$ is uniformly
 bounded by a constant $C_{p_0}$. But then in $x=\sum_{m\bZ} 2^{-m}
 r_m$, $r_m$ is 0 as soon as $2^{-m}>\|x\|$, and
 $z=\sum_{2^{-m}\leq\|x\|}2^{-m} z_m$ is bounded by $2C_{p_0}$.}
    \end{rk}
\begin{rk}{\rm The proof only uses estimates for the norm of $D_{2^k}$ on $L_{p}$. Thus, we may replace in the arguments 
above the space $L_p$ by any symmetric function space $E\subset L_{p_0}+L_{p_1}$ with lower Boyd index $p_E>p$ and upper Boyd index $q_E<p_1$.
This way, we can conclude that $S$ is also
    bounded from $E(\M)$ to $E(\N;\ell_\infty)$. We refer to either of \cite{Boy69, DirkTAMS, DirkJOT, Mon96} for details.
}
    \end{rk}

\begin{rk}{\rm As for the commutative case or in \cite{DirkJOT}, it is also possible to relax to
    $0<p_0< p_1\leq \infty$ with  worst constants. Let us keep the notations introduced in the proof of Theorem \ref{Mar} and write $z_m = \sum_{k\in I} \tilde d_k q_{k,m}$ the sequence obtained in Lemma \ref{basic}. We have $z=\sum_{m\in \bZ} 2^{-m}\sum_{k\in I}\tilde d_k q_{k,m}$ and to estimate its norm, one has to use the $p$-triangular inequality if $p<1$ rather than the Hardy-Littlewood majorization.
    Indeed, we have
    $$
    \norm{z}_p^p \leq 
    \sum_{k\in I ,m\in \bZ} 2^{-mp} \tilde d_k^{\,p} \|q_{k,m}\|_p^p\leq \sum_{k\in I}
  \tilde d_k^{\,p} 2^k   \sum_{m\in \bZ } 2^{-mp}  \tau(r_m)\leq \frac {\|x\|_p^p}{1-2^{-p}}\sum_{k\in I}  \tilde d_k^{\,p} 2^k.$$
 To obtain the last inequality, we used Lemma \ref{decdya} with $\alpha = p$. Then, one has to choose $d_k$ as before so that the series converges. Unfortunately, we have no relevant non commutative applications for the moment.}
\end{rk}

\begin{rk}{\rm As for the commutative case or in \cite{DirkJOT}, we can merely assume that $S_n$ is sub-additive or even that $\mu(S_n(f+g))\preceq \mu(S_n(f))+\mu(S_n(g))$. Indeed, we just use linearity in \eqref{linear} where sub-additivity is enough if one replaces $S_n(r_m)$ by $S_n(2^m r_m)$ in the arguments of Lemma \ref{basic}, there is no need to use homogeneity either. Moreover we only used the hypothesis for projections, so we can change weak type to restricted weak type as in \cite{DirkJOT}.}
\end{rk}

\begin{rk}\label{strongweak}{\rm Assuming the following stronger form of weak type inequality: for any projection $r$
 there is  $a\in L_{p_0,\infty}$ with $\|a\|_{p_0,\infty}^{p_0}\leq \tau (r)$ and 
$0\leq S_n(p)\leq a$ and similarly for $p_0$ with an operator $b$ would remove the square in the constant $C_{\theta,p_0,p_1}$.
Indeed, using the same proof, one can take  $e_k$ to be a spectral projection of $b$
for $k\leq 0$ and a spectral projection of $e_0ae_0$ for $k>0$.
The assumption yields $0\leq S_n(r)\leq 2 \Big((1-e_0) b(1-e_0) + e_0 a e_0\Big)$ without Lemma \ref{basic} being needed so that we simply end up with $0\leq S_n(p)\leq 8 \sum_{k} \min(2^{-k/p_0},2^{-k/p_1})q_k$. }
\end{rk}

\begin{rk}\label{rk:opti}
  {\rm
  The constant in the second part of the theorem is optimal in the following sense: assume that the statement holds with 
  $$
  C'_{\theta,p_0,p_1} \leq CC_0^{1-\theta}C_1^\theta f(\theta,p_0,p_1)
  $$
  for a certain function $f$. Then 
  $$
  \alpha_{\theta,p_0,p_1}^2 = \mathrm{o}(f(\theta,p_0,p_1))
  $$
  when $\theta$ goes to $0$ and when $\theta$ goes to $1$ (unless $q = \infty$, in this case the constant remains bounded when $\theta$ goes to $1$).  A proof,  when $p_0=1$ and $p_1=\infty$, will be proved at the end of the section.
  }
\end{rk}

Our main application is about martingale theory to recover the Doob maximal inequality of \cite{J2}. Assume that $(\M,\tau)$ endowed with an increasing filtration
$(\M_n)_{n\geq 0}$ and associated conditional expectations
$(\E_n)_{n\geq 0}$.

Cuculescu's construction gives that $S=(\E_n)_{n\geq 0}$ is of weak type (1,1),
it is obviously of strong type $(\infty,\infty)$.
Then Theorem \ref{Mar} provides a substitute for  maximal functions
of martingales.

\begin{cor} Let $x\in (L_p(\M) \cap L_1(\M))^+$ for some $1<p\leq \infty$, then there is
  $z\in L_p(\M)^+$ such that
  $$ 0\leq \E_n(x)\leq z \quad \textrm{and} \quad \forall 1<q\leq  p,\; \|z\|_q\leq C_q \|x\|_q$$
  with $C_q=O\big(\mathrm{ln}(q-1)/(q-1)\big)^2$ when $q\to 1$ and $C_q=O(1)$ when $q\to \infty$.
\end{cor}

\begin{proof}
  If $x\in \M^+$, the statement follows directly from Theorem \ref{Mar} when $q<\infty$. For $q=\infty$, this is Remark \ref{casinfini}.

  When $x\in L_p^+$, $1<p<\infty$, this is justified by \eqref{majz}.
\end{proof}

\begin{rk}\label{strongweak2}{\rm Since the optimal behaviour of the constant in the Doob maximal inequality is known to be of order $(p-1)^{-2}$ when $p$ goes to 1, it follows that it not possible to strengthen the weak (1,1)-inequality
    as  in Remark \ref{strongweak}.}
\end{rk}

\subsection{The asymmetric Marcinkiewicz theorem}

In \cite{J2}, asymmetric maximal inequalities were considered, they can be deduced in the same way.

For convenience, we recall the definition of asymmetric maximal function spaces. For $0<p\leq\infty$ and $0<\gamma<1$, $L_p(\N;\ell_\infty^\gamma)$ consists of all
sequences of operators $x=(x_n)_{n\geq0}$ for which there exist operators 
$a,b\in L_p(\N)^+$ and contractions $(u_n)_{n\geq0}$ in $\N$ such that
\begin{equation*}
    \forall n\geq1,\qquad x_n =a^\gamma u_nb^{1-\gamma}.
    \end{equation*}
The associated norm is $\|x\|_{L_p(\N;\ell_\infty)}=\inf \|a\|_{p} \|b\|_{p}$, 
where the infimum runs over all decompositions as above. Then $L_p(\N;\ell_\infty^\gamma)$ is a Banach space for $p\geq\max\{2\gamma,2(1-\gamma)\}$ and a quasi-Banach otherwise. Of course, we recover $L_p(\N;\ell_\infty)$ when $\gamma=\frac 12$. The limit cases $\gamma=0,1$ correspond to the column and row maximal function spaces $L_p(\N;\ell_\infty^c)$ and $L_p(\N;\ell_\infty^c)$ that are recalled at the beginning of section \ref{sec:beyondpos}.

For a sequence of maps $S=(S_n)$ as before, we say that it is of strong $\gamma$-asymmetric type $(p,p)$ if $S$ is bounded from $L_p(\M)$ to $L_p(\N;\ell_\infty^\gamma)$. We easily get
\begin{thm}\label{thm:asym}
  Let $1\leq p_0< p_1\leq \infty$ and
$S=(S_n)$ be a sequence of linear positive maps from 
$L_{p_0}(\M)+ L_{p_1}(\M)$ to $L_0(\N)$. Assume that $S$ if of weak type 
$(p_0,p_0)$  and $(p_1,p_1)$ with constants $C_0$ and $C_1$.

Then for any $0<\theta,\gamma <1$,
 then $S$ is of strong $\gamma$-asymmetric $(p,p)$ type
 where $\frac 1 {p}= \frac{1-\theta}{p_0}+\frac \theta{p_1}$
 if $p>\max\{2\gamma,2(1-\gamma)\}$.

Moreover under these conditions, there is a constant
$C_{\theta,\gamma,p_0,p_1}$ such that for any $x\in L_{p_0}(\M)\cap
L_{p_1}(\M)$, there exist $a_\gamma, b_\gamma \in L_{p}^+$
 and contractions $u_n\in \M$ (all independent of $\theta$) with
$$S_n(x)= a_\gamma^\gamma u_n b_\gamma^{1-\gamma},\qquad  \|a_\gamma\|_{p},\|b_\gamma\|_{p}\leq  C_{\theta,\gamma,p_0,p_1} \|x\|_{p}.$$
\end{thm}

\begin{proof}
  This is just a variation on the previous argument. We use the notation from the proof of Theorem \ref{Mar}.

 We assume $x\in L_{p_0}(\M)\cap L_{p_1}(\M)$ is positive. Writing it
 $x=\sum_{m\in \bZ} 2^{-m} r_m$. We fix $d_k={|k|(\ln|k|)^2+1}$ to construct the elements $z_m$ in Lemma \ref{basic}.

 We have that there exist contractions $c_{m,n}\in \M$ such that
 $S_n(r_m)=z_m^{1/2} c_{m,n} z_m^{1/2}$. We use Lemma \ref{fact}
to get contractions $v_{m,n}$ such that
$$S_n(x)= \sum_{m\in \bZ} 2^{-\gamma m} z_m^{1/2} c_{m,n}
2^{-(1-\gamma) m} z_m^{1/2} =\Big(\sum_{m\in \bZ} 2^{-2\gamma m}
z_m\Big)^{1/2} v_{m,n} \Big(\sum_{m\in \bZ} 2^{-2(1-\gamma) m}
z_m\Big)^{1/2}.$$ We set $a_\gamma= \Big(\sum_{m\in \bZ} 2^{-2\gamma
  m} z_m\Big)^{1/(2\gamma)}$ and $b_\gamma=\Big(\sum_{m\in \bZ}
2^{-2(1-\gamma) m} z_m\Big)^{1/(2(1-\gamma))}$. We now justify that
they are in $L_{p}$ that also legitimates the use of  Lemma \ref{fact}.
Thanks to Lemma \ref{decdya}
$$\mu(a_\gamma)^{2\gamma}\preceq \sum_{m\in \bZ} 2^{-2m\gamma} \mu(z_m)\preceq
 \sum_{m\in \bZ}\sum_{k\in \bZ} 2^{-2m\gamma}\tilde d_k D_{2^k} (\mu(r_m))\leq
C_\gamma \sum_{k\in \bZ}\tilde d_k D_{2^k} (\mu(x^{2\gamma})). $$
And since $p/(2\gamma)\geq 1$, we can use the triangular inequality to get
$$\|a_\gamma\|_{p}=\|a_\gamma^{2\gamma}\|_{p/(2\gamma)}^{1/(2\gamma)}\leq
\Big(C_\gamma\sum_{k\in \bZ} \tilde d_k 2^{k/p}\|x^{2\gamma}\|_{p/(2\gamma)}\Big)^{1/(2\gamma)}\leq C_{\theta,\gamma,p_0,p_1} \|x\|_{p}.$$
One deals with $b_\gamma=a_{1-\gamma}$ in the same way. 
\end{proof}

The condition $p>\max\{2\gamma,2(1-\gamma)\}$ can not be removed in general, we provide an easy example for $p_0=1$ and $p_1=\infty$.

We choose $\M=\bC$ and $\N=\mathbb B(\ell_2)$ with their natural traces. Let
$S_n(\lambda)=\lambda T_n$ with $T_n=e_{1,1}+1/\sqrt{n} (e_{1,n}+e_{n,1})+1/n e_{n,n}\geq0$.

It is clear that $(T_n)_{n\geq 1}$ is bounded in $\mathbb B(\ell_2)$, thus $S$ is of strong $(\infty,\infty)$ type.

For any $t>0$, the projection $r=\sum_{k>1/t} e_{k,k}$ satisfies
$\tau(1-r)\leq  1/t$ and $\|r T_n r\|_\infty\leq 2 t$. Thus $S$ is also of weak (1,1) type.

\begin{prop} \label{prop:noasym}
For $\theta>1/2$, $S$ is not of strong $\theta$-asymmetric $(2\theta,2\theta)$ type.
  \end{prop}
\begin{proof}
  Otherwise, there would exist $a,b\in S_{2\theta}^+$ and contractions $u_n$ such that $T_n=a^\theta u_n b^{1-\theta}$.
  Denoting by $(\delta_i)$ the canonical basis of $\ell_2$, we would have
  $$\langle a^\theta \delta_n, u_nb^{1-\theta}\delta_1\rangle =1/\sqrt{n}\leq \| a^\theta \delta_n\|.\|b^{1-\theta}\delta_1\|.$$
  Thus $\| a^\theta \delta_n\|\geq C/\sqrt{n}$ yielding  $\|a^\theta\|_2^2\geq \sum_{n\geq 1 }C/n=\infty$ and $a_\theta\notin S_{2\theta}$, a contradiction.
\end{proof}

We conclude by noting that Theorem \ref{thm:asym} also holds for $\gamma=0,1$ in the case $p_0\geq1$ and $p_1=\infty$ if $p>2p_0$. Indeed, if $x\in L_p^+$, then $S_n(x)^2\leq C_1 S_n(x^2)$ (because $p_1=\infty$) and by Theorem \ref{Mar}, $S_n(x)^2\leq a^2$ for some $a\in L_{p}$. Thus we may conclude that 
$S_n(x)=u_na=au_n^*$ for some contractions $u_n$ using the polar decomposition. This is similar to \cite{J2} section 5. It is not possible to go down to $p=2p_0$ with the previous counterexample when $p_0=1$ (or a variation if $p_0>1$).

\subsection{Constructing counterexamples}

We conclude this section by presenting a proof of Remark \ref{rk:opti} and a lemma that will be used to construct other counterexamples in the next section. Recall than Lemma \ref{majdiag} introduced a way to construct a supremum for families of operators for which a diagonal is controlled. The following lemma shows that this construction is, in a certain sense, optimal. 

\begin{lemma}\label{lem:counter}
  Let $\H$ be a Hilbert space. Let $(p_i)_{i\leq N}$ be a finite family of orthogonal projections in $\B(\H)$ and $(\alpha_i)_{i\leq N}$ positive real numbers. Let $a \in B(\H)^+$ be such that  $a^2 \geq b^*b$ for any $b$ of the form 
  $$
  b = \sum_{i\leq N} \alpha_ic_ip_i, \quad \norm{c_i} \leq 1.
  $$
  Set $p = \sum_{i\leq N} p_i$. Then there exists an invertible contraction $C$ and a sequence of positive real numbers $(\lambda_i)_{i\leq N}$ such that 
  $$
  ap = \dfrac12C^{-1}(\sum_{i\leq N} \lambda_i^{-1/2}\alpha_ip_i) \quad \text{and} \quad \sum_{i\leq N} \lambda_i = 1.
  $$
\end{lemma}

\begin{proof}
  For $i\leq N$, set $\H_i = p_i\H$. Without loss of generality, we can assume that $\H = \bigoplus_{i\leq N} \H_i$. Set $\tilde a=a \big(\sum_{i\leq N} \alpha_i^{-1}p_i\big)$, so that $a^2 \geq b^*b$ for $b = \sum_{i\leq N} c_ip_i$ where $c_i$ are contractions.
  
 Let $\xi = (\xi_i)_{i\leq N} \in \H$, by choosing contractions $c_i$ such that 
  $$
  \Big\|{\sum_{i\leq N} c_i\xi_i}\Big\|_{\H} = \sum_{i\leq N} \norm{\xi_i}_{\H_i},
  $$
   we obtain
  $$
  \norm{a\xi}_\H \geq \sum_{i\leq N} \norm{\xi_i}_{\H_i}.
  $$
  Consequently $\tilde a^2\geq 1$ so $\tilde a$ is invertible and its inverse $\tilde a^{-1}$ can be regarded as contraction from $\H$ to $\ell_1((\H_i)_{i\leq N})$. Denote by $c$ the anti-adjoint of $\tilde a^{-1}$, it is contractive from  $\ell_\infty((\H_i)_{i\leq N})$ to $\H$. Since $\H_i$ can sits as a $1$-complemented subspace in $\B(\H_i)$, $c$ can be extended to a contraction from $A = \ell_\infty((\B(H_i))_{i\leq N})$, which is a $C^{*}$-algebra, to $\H$. By the little Grothendieck theorem \cite[Theorem 9.4]{Pi86}, there exists a state $\varphi$ on $A$ such that for any $x = (x_i)_{i\leq N} \in A$, 
 $$
  \norm{c(x)}_{\H}^2 \leq 2\varphi(xx^* + x^*x).
  $$
  The state $\varphi$ can be decomposed $\varphi(x) = \sum_{i\leq N} \lambda_i\varphi_i(x_i)$ where $\varphi_i$'s are  states on $\B(\H_i)$ and  $\lambda_i$'s are positive real numbers such that $\sum_{i\leq N} \lambda_i = 1$. Then 
  $$
  \norm{c(x)}_{\H}^2 \leq 2\sum_{i\leq N} \lambda_i\varphi_i(x_ix_i^* + x_i^*x_i) \leq 4\sum_{i\leq N} \lambda_i\norm{x_i}_{B(\H_i)}^2.
  $$
  Consider now $c$ as a bounded operator from $\H \subset A$ to $\H$. The previous inequality shows that $c^*c \leq 4d^2$ where $d$ is the diagonal operator $d = \sum_{i\leq N} \sqrt{\lambda_i} p_i$. So $c$ can be factorized as $c = 2Cd$ where $C$ is a contraction. Hence, $\tilde a^{-1}$ admits a factorization of the form $\tilde a^{-1} = 2dC'$, which concludes the proof going back to $a$. 
\end{proof}

\begin{proof}[Proof of Remark \ref{rk:opti}]
Set $\M = \bC$ and $\N = \B(\ell_2)$. For any $n\geq 0$, let 
$$
q_n = \sum_{i=2^n-1}^{2^{n+1}-2} e_{i,i} \in \B(\ell_2) \quad \textrm{and} \quad Q_n= \sum_{n= 0}^N q_n. 
$$
Let $X \subset \B(\ell_2)^+$ be the set of operators $x$ of the form
$$
x = b^*b \quad \text{with} \quad b = \sum_{i=0}^\infty 2^{-i/2}c_iq_i, \quad \norm{c_i} \leq 1.
$$
Consider the family of operators $S = (S_x)_{x\in X}:\bC \to \B(\ell_2)$ defined by 
$$
S_x : t \in \bC \mapsto tx.
$$
By basic inequalities $\|Q_n^\bot x Q_n^\bot\| \leq 2^{-n}$ and ${\rm tr}\, Q_n= 2^n$. Hence, one $S$ is of weak type $(1,1)$ and of strong type $(\infty,\infty)$. Let $a \in \B(\ell_2)$ be such that $a \geq S_x(1)$ for any $x\in X$, meaning $a\geq x$ for any $x\in X$. We want an estimate of $\|a\|_p$ as $p\to 1$.

By Lemma \ref{lem:counter}, for any $N>0$, there exists a sequence of positive reals $(\lambda_{i,N})_{0\leq i\leq N}$ such that
$$
4Q_NaQ_N \geq \sum_{i=0}^N 2^{-i}\lambda_{i,N}^{-1} q_i \quad \text{and} \quad \sum_{i=0}^N \lambda_{i,N}= 1.
$$
Fix an ultrafilter $\mathfrak U$ on $\mathbb N$ and let $\lambda_i=\lim_{N, \mathfrak U} \lambda_{i,N}$. We have $\sum_{i=0}^\infty \lambda_{i}\leq 1$. Since $\lambda_{i,N}^{-1}\leq 2^{i+2}\|a\|$, we can also conclude that $\lambda_{i}>0$ and that  for all $N$, 
$4Q_NaQ_N \geq \sum_{i=0}^N 2^{-i}\lambda_{i}^{-1} q_i$.

Clearly, we have that $4^p\|a\|_p^p\geq \sum_{i\geq 0} \lambda_{i}^{-p}2^{i(1-p)}$. Thus using basic inequalities and the Cauchy-Schwarz inequality
$$\|a\|_p^p\gtrsim \sum_{i= 1/2(p-1)}^{1/(p-1)} {\lambda_i}^{-p}\gtrsim \frac 1{(p-1)^2}\Big(\sum_{i\geq 1/2(p-1)}\lambda_i^p\Big)^{-1}\geq  \frac {R_p^{-1}}{(p-1)^2},$$
where $R_p=\sum_{i\geq 1/2(p-1)}\lambda_i$. This leads to a constant of order strictly higher than $(p-1)^{-2}$ as $p$ goes to $1$.
\end{proof}

One can also prove Remark \ref{rk:opti} using averaging techniques by showing that $a$ can be chosen of the form $\sum d_i q_i$ ($d_i\in \mathbb R^+$) as long as $p\geq 1$.

\section{Beyond positivity: $\Lambda$-spaces}\label{sec:beyondpos}

The positivity assumption  for $S$ can not be removed in Theorem \ref{Mar}. Indeed, consider   $\M=\bC$ and $\N=\mathbb B(\ell_2)$ with their natural traces as above and set $S_n(\lambda)=\lambda T_n$ with $T_n=e_{n,1}+ e_{1,n}$. The map $S=(S_n)$ is clearly of strong type $(\infty, \infty)$ and weak type (1,1). If $S$ was of strong type $(p,p)$, this would imply that there is $A\in S_p^+$ such that $-A \leq T_n \leq A$. This would force that $A_{1,1}A_{n,n}\geq 1$ for all $n\geq 1$, this is impossible as we must also have $A_{n,n}\to 0$.

It is however possible to get some results but with weaker factorizations using the row and column
$\ell_\infty$-valued Lorentz spaces and some variations. Our reference is \cite{XW}.

Let $1\leq p,q \leq \infty$, a sequence $x=(x_n)$ of elements in
$L_{p,q}(\N)$ belongs to $L_{p,q}(\N;\ell_\infty^c)$ (or simply
$L_{p,q}(\ell_\infty^c)$) if there is some $a\in L_{p,q}(\N)^+$ and
contractions $u_n\in \N$ such that $x_n=u_na$. This is equivalent to say that
$x_n^*x_n\leq a^2$ for all $n\geq 0$. The infimum of
$\|a\|_{p,q}$ over all possible $a$ defines the quantity $\|x\|_{L_{p,q}(\ell_\infty^c)}$. We obtain a Banach space if $p>2$ or if $p=2$ and $q\geq 2$. 

The row version  $L_{p,q}(\N;\ell_\infty^r)$ is obtained by taking adjoints.

The couples $(L_q(\N;\ell_\infty^c), L_r(\N;\ell_\infty^c))$ are
compatible in the sense of interpolation theory (we assume $r>q$). One of the result in \cite{Mu} and generalized in \cite{XW} is that they behave well with
respect to the complex interpolation method (in the Banach spaces range).

\medskip

We will also need a weaker version of non commutative maximal
inequalities. The construction naturally extends the definition of
weak type $(p,p)$ for maximal operators and still coincides with
standard maximal inequalities in the commutative case. It has the
advantage of retaining its relation with almost uniform and bilateral
almost uniform convergence so it may be used as a substitute for the
strong version of maximal inequalities when studying pointwise
convergence questions occurring in ergodic theory and Fourier
analysis.

Our starting point is the following notions of non-increasing rearrangement for a sequence $(x_n)_{n\geq 0}$ that should be thought of as a non commutative analog of $\mu(\sup_{n\geq0} \md{x_n})$.

\begin{defi}
  Given a sequence $x = (x_n)_{n\geq 0}$ in $L_0(\N)$, we define three non-increasing functions 
  $\mu(x), \mu_c(x), \mu_r(x): \bR^{+*}\to \bR^+$ as, for $t>0$
  $$
  \mu(x,t)= \inf_{\tau(1-e)\leq t}  \sup_{n\geq 0} \| ex_ne\|,\quad 
 \mu_c(x,t)= \inf_{\tau(1-e)\leq t}  \sup_{n\geq 0} \| x_ne\|,
 \quad \mu_r((x_n)) = \mu_c((x_n^*))
 $$ where the infimun runs over all projections $e\in \N$.
\end{defi}
When $x$ is a constant self-adjoint sequence, $x_n= a = a^*$ for all $n\geq0$, we recover $\mu(x)=\mu(a)$. We
could have used pairs of projections $e,f$ with
$\tau(1-e),\tau(1-f)\leq t$ and $\| fx_ne\|$ to fully recover $\mu(a)$
for general $a$ but this won't make a significant difference. Note
also that the following fundamental inequality is still verified for
any $s,t>0$ and $x = (x_n), y = (y_n) \in L_0(\N)^{\mathbb{N}}$
  \begin{equation}\label{eq:mulambda}
      \mu_\sharp(x+y,t+s) \leq \mu_\sharp(x,t) + \mu_\sharp(y,s),
    \end{equation}
    with $\sharp=c,\,r,\,\emptyset$.  This also holds with $s=0$
    replacing $\mu_\sharp(y)(s)$ by $\|y\|$.

  The definition is motivated by weak type maximal inequalities and made so that $\|(x_n)\|_{\Lambda_{p,\infty}(\N;\ell_\infty)}=\sup_{t>0} t^{1/p} \mu((x_n))(t)$. Indeed, one just need to set $\lambda= \frac C{t^{1/p}}$ in \eqref{normefaible}. Thus
  $$ \Lambda_{p,\infty}(\N;\ell_\infty)= \{ (x_n)\in L_p(\N)^\bN  \;|\;  \|(x_n)\|_{\Lambda_{p,\infty}}= \sup_{t>0} t^{1/p} \mu((x_n))(t)<\infty\}.$$
  
  It is natural to extend the definition to $L_p$-spaces or Lorentz spaces. 
  \begin{defi}
    For $p,q>0$, we define for $\sharp=c,\,r,\,\emptyset$ $$\Lambda_{p,q}(\N;\ell_\infty^\sharp)= \{ (x_n)\in L_p(\N)^\bN  \;|\;  \|(x_n)\|_{\Lambda_{p,q}^\sharp}= \|\mu_\sharp((x_n))\|_{p,q}<\infty\}.$$
  \end{defi}
  
  As usual, we will write $\Lambda_p$ instead of $\Lambda_{p,p}$. We may also write $\Lambda_{p,q}^\sharp$ instead of $\Lambda_{p,q}(\N;\ell_\infty^\sharp)$ to lighten notations. Of course $\Lambda_\infty^\sharp=L_\infty(\N;\ell_\infty^\sharp)=\ell_\infty(\N)$. Let us collect a few simple properties of these spaces.  
  
  \begin{prop}\label{quasi}
  Let $p,q  \in (0,\infty]$.
   \begin{enumerate}
       \item $\Lambda_{p,q}(\N;\ell_\infty^\sharp)$ is a  quasi-Banach space
       \item for any sequence $x \in L_p(\M;\ell_\infty^\sharp)$, 
 $$
 \norm{x}_{\Lambda_p(\N;\ell_\infty^\sharp)} \leq 2^{1/p}\norm{x}_{L_p(\M;\ell_\infty^\sharp)}
 $$
   \end{enumerate}
  \end{prop}
  
  \begin{proof}
Point (1) is clear using \eqref{eq:mulambda}. Let us prove (2). Let $X = (x_n)_{n\geq0}$ be a sequence in $L_p(\N;\ell_\infty)$. Let $a, b \in L_{p/2}(\N)$ and $(y_n)_{n\geq 0}$ in $\N$ a sequence of contraction such that for any $n\geq 0$, $x_n = ay_nb$. Let $t > 0$. We can find a projection $e\in \N$ such that $\tau(e_1) \leq t$ and $\norm{ea}_\infty \leq \mu(a,t)$. Similarly, we choose a projection $e_2$ such that $\tau(e_2) \leq t$ and $\norm{be_2}_\infty \leq \mu(b,t)$. Set $e = e_1 \vee e_2$. It is clear that $\tau(e_2) \leq 2t$ and for any $n\geq 0$, $\norm{ex_ne}_\infty \leq \mu(a,t)\mu(b,t)$. Therefore 
$$
\mu(X,2t) \leq \mu(a,t)\mu(b,t)
$$
which implies the desired inequality by integrating over $t$ and H\"older's inequality. 

When $\sharp=c,r$, the proof is simpler and there is actually no factor $2^{1/p}$.
\end{proof}

Controlling the $\Lambda_p$ norm of a sequence $(x_n)_{n\geq 0}$ is  much weaker than controlling its standard maximal norm. In particular, it does not allow to exhibit an element in $L_p(\N)$ that would play the role of $\sup_{n\geq 0} x_n$, even for positive sequences.

Recall that we use the classical notation from \cite{BL} concerning interpolation theory.  Since every $\Lambda_{p,q}^\sharp$  can be continuously embedded in the topological vector space $L_0(\N)^{\mathbb{N}}$, these spaces are all compatible in the sense of interpolation. The following proposition asserts that the $\Lambda_{p,q}^\sharp$ form a (real) interpolation scale.
\begin{prop}
  For $0<p<\infty$, $t>0$ and any $x=(x_n)\in \Lambda_p^\sharp+ \Lambda_\infty^\sharp$:
  $$K(t, x, \Lambda_p^\sharp,\Lambda_\infty^\sharp) \simeq_p K(t, \mu_\sharp(x), L_p, L_\infty).$$
  Consequently $(\Lambda_{p_0,q_0}^\sharp,\Lambda_{p_1,q_1}^\sharp)_{\theta,q}\simeq\Lambda_{p,q}^\sharp$ for $0<\theta<1$ with $\frac 1 p=\frac{1-\theta}{p_0} +\frac \theta{p_1}$, $0<p_0<p_1\leq\infty$ and $0<q\leq\infty$.
\end{prop}
\begin{proof}
  We detail the proof only in the case $\sharp=\emptyset$, others are similar.
  
   First, fix $t>0$, $\eps>0$ and choose a decomposition $x=a+b$ such
   that
   $\|a\|_{\Lambda_p}+t \|b\|_{\Lambda_\infty}\leq (1+\eps) K(t, x,
    \Lambda_p, \Lambda_\infty)$. Since
   $\|b\|_{\Lambda_\infty}= \|b\|_\infty$, we have that
   $\mu(x)(t)\leq \mu(a)(t)+\|b\|_\infty$. Hence, we can decompose
   $\mu(x)=\alpha+\beta$ with $\alpha\leq \mu(a)$ and
   $\beta\leq \|b\|_\infty$. It follows that
   $K(t, \mu(x), L_p, L_\infty)\leq \| a\|_{\Lambda_p}+
   t\|b\|_{\Lambda_\infty}$. Thus
   $K(t, \mu(x), L_p, L_\infty)\leq K_t(t,x, \Lambda_p, \Lambda_\infty)$.

   For the opposite, we use the Holmstedt formula $K(t, \mu(x), L_p, L_\infty)\simeq_p \Big(\int_0^{t^p} \mu(x)^p(s)ds\Big) ^{1/p}$ for $t>0$. Fix $t>0$, we can find a projection $e\in \N$ such that $\tau(1-e)\leq t^p$ and $\|ex_ne\|\leq 2 \mu(x)(t^p)$. We decompose $x=exe+\big((1-e)x+ex(1-e)\big)=b+a$. If $s>t^p$ then
   clearly $\mu(a)(s)=0$ and otherwise, $\mu(a)(s)\leq \mu(x)(s)+2\mu(x)(t^p)$.
   We get
   $$\|\mu(a)\|_p^p \lesssim_p \int_0^{t^p} \mu(x)^p(s)ds + t^p  \mu(x)^p(t^p)\leq 2 \int_0^{t^p} \mu(x)^p(s)ds.$$
   Thus $ K_t(t,x, \Lambda_p, \Lambda_\infty)\leq \|\mu(a)\|_p + t  2 \mu(x)(t^p)\lesssim_p \Big(\int_0^{t^p} \mu(x)^p(s)ds\Big) ^{1/p}$. This concludes the proof for the $K$-functional.

   The last statement follows from the estimate for the $K$-functional and the
   corresponding result for commutative $L_p$ and the reiteration principle.
 \end{proof}
 The following version of the Marcinkiewicz theorem is now clear:
\begin{prop}
 Assume that $S$ is of weak type $(p_0,p_0)$ and $(p_1,p_1)$ with $p_0, p_1\geq 1$. Then $S: L_p(\M)\to \Lambda_p(\N;\ell_\infty)$ is bounded. 
\end{prop}

The spaces $\Lambda_{p,q}^c$ are connected to $L_{p,q}(\ell_\infty^c)$ via real interpolation. To make this precise, we start with an effective  characterization
of $\Lambda_{p}^\sharp$.

To improve clarity, we use capital letter for sequences. Thus
recall that if $X=(x_n)$, for an element $n\in \N$, we set $Xq=(x_nq)$
and similarly on the left. We will consider the weighted Lorentz spaces $\ell_{p,q}^\omega$ associated to the measure on $\bZ$ given by $\omega(\{n\})= 2^n$.

\begin{lemma}\label{decomp} We have the following for $0<p<\infty$ and $0<q\leq \infty$:
  
  \begin{enumerate}
  \item\label{disc}  A sequence $X\in \Lambda_{p,q}^c$ if and only if there exist a sequence of disjoint projections $q_k\in \N$ with $\tau(q_k)\leq 2^{k}$, $(a_k)\in \ell_{p,q}^{\omega,+}$  and contractions $U_k=(u_{n,k})\in \ell_\infty(\N)$ for $k\in \mathbb Z$ such that $X=\sum_{k\in \bZ} a_k U_kq_k$.   Moreover $\|X\|_{{\Lambda_{p,q}^c}}\simeq_{p,q} \|(a_k)\|_{p,q,\omega}$.
    \item\label{sum}  We have $\Lambda_{p,q}=\Lambda_{p,q}^c+\Lambda_{p,q}^r$ with equivalent semi-norms.
  \end{enumerate}
\end{lemma}
 \begin{proof}
   We start with the left to right implication in (\ref{disc}). We may assume that
   $\|X\|_{\Lambda_{p,q}^c}=1$. From the definition of $\mu_c$, we may
   find projections $f_k$ such that $\tau(1-f_k)\leq 2^k$ and
   $\|Xf_k\|\leq \mu_c (X)(2^{k})+ 2^{-k^2}=a_k$.
   As we did in Lemma \ref{basic}, we may replace $(f_k)$ by a sequence of smaller decreasing projections $(e_k)$ such that $\tau(1-e_k) \leq 2^{k+1}$. Since
   $\sum_{k\in \bZ}(\mu_c(2^k) + 2^{-k^2}) 1_{(2^{k},2^{k+1})}$ on $\bR$  and $(a_k)$ on $(\bZ,\omega)$ have the same distribution up to a dilation by 2, we deduce that $\|(a_k)\|_{p,q,\omega}\leq C_p$, where the constant $C_p$ depends only on $p$ and goes to infinity only when $p$ goes to $0$.
   Set $q_k=e_k-e_{k-1}$, they are disjoint. We clearly have the decomposition
   with $U_k=a_k^{-1} Xq_k$. The convergence of the series holds in $L_0(\N)^{\mathbb N}$ for instance. 

   For the other implication, we prove it first when $q=p$ and $p> 1$. Assume that $X=\sum_{k\in \bZ} a_k U_kq_k$. Clearly $$\mu_c(X,2^l)\leq \big\|\sum_{k\geq l} a_kU_k q_k\big\| \leq \sum_{k\geq l} a_k\leq \Big(\sum_{k\geq l} 2^{ k}a_k^p\Big)^{1/p}
   \Big(\sum_{k\geq l} 2^{-k/(p-1)}\Big)^{(p-1)/p},$$
  thanks to the H\"older inequality. We get 
   $$
   \|\mu_c(X)\|_p^p
   \leq \sum_{l\in \bZ} 2^l\mu_c(X,2^l)^p
   \leq \frac{1}{(1-2^{-1/(p-1)})^{p-1}} \sum_{l\in \bZ}\sum_{k\geq l} 2^{k}a_k^p
   \leq K_p\sum_{k\in \bZ}2^{k }a_k^p. 
   $$
 
Note that $K_p^{1/p}$ remains bounded when $p\to 1$ but goes to $\infty$ with $p$.

If $p\leq 1$, one just need to use $ (\sum_{k\geq l} a_k)^p\leq  \sum_{k\geq l} a_k^p$ to get the same estimate with $C_p=1$.

The Lorentz case follows easily by interpolation using the linear map
$(a_k)\mapsto \sum_{k\in \bZ} a_k U_kq_k$.

\medskip

To deal with (\ref{sum}), first the bounded inclusion
$\Lambda_{p,q}^c+\Lambda_{p,q}^r\subset \Lambda_{p,q}$ is clear as
$\mu\leq \mu_c,\mu_r$. We use (\ref{disc}) for the reverse. Let
$X\in \Lambda_{p,q}$ with norm 1.  Just as we did, we can find a
sequence of decreasing projections $f_k$ such that
$\tau(1-f_k)\leq 2^k$ and
$\|f_kXf_k\|\leq \mu (X)(2^{k})+ 2^{-k^2}=a_k$. Set $q_k=f_k-f_{k-1}$
and as $X=\sum_{k\in\bZ} f_kXq_k+f_{k+1}Xq_k$, decompose
$$ X =\sum _{k\in \bZ} a_k   (a_k^{-1}f_k Xq_k) q_k +  a_k   q_k(a_k^{-1}q_k Xf_{k+1})= \sum  _{k\in \bZ} a_k (U_k q_k + q_k V_k).$$ 
Clearly $U_k$, $V_k$ are contractions. By the first point $\sum  _{k\in \bZ} a_k U_k q_k\in \Lambda_{p,q}^c$ with norm less than $C_{p,q}$ and similarly for  $\sum  _{k\in \bZ} a_k q_k V_k \in \Lambda_{p,q}^r$.

\end{proof}

\begin{rk}\label{decompalpha}
  {\rm We choose to use $(2^k)$ for simplicity but the above lemma works more generally for any geometric sequence $(2^{\alpha k})$ with $\alpha>0$.
    Namely,  $X\in \Lambda_{p,q}^c$ iff it can be decomposed $X=\sum_{k\in \bZ} a_k U_kq_k$ where $(q_k)$ is a sequence of disjoint projections with $\tau(q_k)\leq 2^{\alpha k}$, $(a_k)\in \ell_{p,q}^{\omega^\alpha,+}$  and $U_k=(u_{n,k})\in \ell_\infty(\N)$ are contractions for $k\in \mathbb Z$.
    Moreover $\|X\|_{{\Lambda_{p,q}^c}}\simeq_{p,q,\alpha} \|(a_k)\|_{p,q,\omega^\alpha}$.}
\end{rk}
\begin{cor}\label{interell}
  Let $0<p<\infty$, $0<q\leq \infty$ and $0<\theta<1$, set $\frac 1 {p_\theta}= \frac {1-\theta} p$,  then with equivalent semi-norms:
  $$ \big(L_p(\ell_\infty^c), L_\infty(\ell_\infty^c)\big)_{\theta,q}= \Lambda_{p_\theta,q}(\ell_\infty^c).$$
 \end{cor}
 \begin{proof}
   As we saw in Proposition \ref{quasi}, the Markov inequality gives contractive inclusions
   $L_p(\ell_\infty^c)\subset \Lambda_p(\ell_\infty^c)$  for
   $p>0$. Thus interpolation immediatly gives a bounded
   inclusion
   $$\big(L_p(\ell_\infty^c),
   L_\infty(\ell_\infty^c)\big)_{\theta,q}\subset\big(\Lambda_p(\ell_\infty^c),
   \Lambda_\infty(\ell_\infty^c)\big)_{\theta,q}\simeq
   \Lambda_{p_\theta,q}(\ell_\infty^c).$$ For the reverse, first take
   $X\in \Lambda_{p_\theta}(\ell_\infty^c)$ and consider the
   decomposition from Lemma \ref{decomp} \eqref{disc} and Remark \ref{decomp} with $\alpha= p$ . We have written
   $X=\sum a_k U_kq_k$ in $L_0(\N)^{\mathbb N}$ with
   $$J(2^k, X_k, L_p(\ell_\infty^c), L_\infty(\ell_\infty^c))=\max \{ \| a_k U_kq_k\|_{L_p(\ell_\infty^c)}, 2^k \| a_k U_kq_k\|_{L_\infty(\ell_\infty^c)}\} \leq a_k 2^{k}.$$
   Moreover $\big(\sum_{k\in\bZ} 2^{kp} a_k^{p_\theta}\big)^{1/p_\theta} \lesssim\|X\|_{\Lambda_{p_\theta}^c}$. Since
   $$\sum_{k\in \bZ} \big(2^{-k\theta} J(2^k, X_k, L_p(\ell_\infty^c), L_\infty(\ell_\infty^c)\big)^{p_\theta}\leq \sum_{k\in \bZ}2^{kp} a_k^{p_\theta},$$
   we conclude that $X\in \big(L_p(\ell_\infty^c),
   L_\infty(\ell_\infty^c)\big)_{\theta,p_\theta}$ by the equivalence between the $J$- and $K$-methods (for quasi-normed spaces). Thus $\big(L_p(\ell_\infty^c),
   L_\infty(\ell_\infty^c)\big)_{\theta,p_\theta}\simeq \Lambda_{p_\theta}(\ell_\infty^c)$. The general statement follows by the reiteration theorem.
\end{proof}

More generally, choosing the correct weight $\omega^\alpha$, the same proof gives that with $0<p_0<p_1$ and $0<\theta<1$ with $\frac 1 p=\frac {1-\theta}{p_0}+\frac 1 {p_1}$ and $0<q\leq \infty$: $$\big(L_{p_0}(\ell_\infty^c), L_{p_1}(\ell_\infty^c)\big)_{\theta,q}= \Lambda_{p_\theta,q}(\ell_\infty^c).$$

\begin{cor}\label{interp}
  For $2<p<\infty$ and $1\leq q\leq \infty$, $\Lambda_{p,q}(\ell_\infty^\sharp)$ has an equivalent norm. 
\end{cor}
\begin{proof}
  For $\sharp=c$, this is clear from Corollary \ref{interell} as $L_2(\ell_\infty^c)$ is a Banach space. The case $\sharp=r$ is obtained by taking adjoints. The remaining case then follows using Lemma \ref{decomp} \eqref{sum}.
\end{proof}

At this point, it is worth to justify that $\Lambda_p(\ell_\infty^c)\neq L_p(\ell_\infty^c)$:

\begin{prop}\label{Ll}
  Set $\N=\mathbb B(\ell_2)$. Let $0<p<\infty$, and $q>2$. The formal
  identity map on $L_p(\ell_\infty^c)$ is
  not bounded from  $\Lambda_{p,q}(\ell_\infty^c)$ to
  $L_{p,\infty}(\ell_\infty^c)$ nor from $\Lambda_{p}(\ell_\infty^c)$ to
  $L_{p}(\ell_\infty^c)$.  
 \end{prop}

\begin{proof}
 Set $N>0$ and let $p_i=\sum_{k=2^{ i -1}}^{2^{ i}-1} e_{k,k}$ for $i=1,...,N$. Choose some
 families of contractions $(u_{n,i})_{n\geq0}$, $i\leq N$ such that $\{(u_{n,1},...u_{n,N});n\geq 0\}$ is dense in $B^N$ where $B$ is the unit ball of compact operators. Set $X_i=2^{-i/p }(u_{n,i}p_i)_{n\geq 0}\in L_p(\ell_\infty^c)\cap L_\infty(\ell_\infty^c)$.

By Lemma \ref{decomp}, the norm  of $X=\sum_{i=1}^N X_i$ in $\Lambda_{p,q}^c$ is controlled by $C_{p,q} N^{1/q}$.
Similarly, we also have $\|X\|_{\Lambda_{p}^c}\lesssim_p N^{1/p}$.

Next we try to estimate the norm of $X$ in $L_{p,\infty}(\ell_\infty^c)$ or in $L_p(\ell_\infty)$.

By definition, we can find $a\in S_{p,\infty}^+$ (or $a\in S_p$) such that $\|a\|_{p,\infty}\leq 2
\| X\|_{L_{p,\infty}(\ell_\infty^c)}$ (or $\|a\|_p\leq 2\| X\|_{L_{p}(\ell_\infty^c)}$) and 
$$ \Big|\sum_{i=1}^N 2^{-i /p} u_{n,i} p_i\Big|^2=
\Big(\sum_{i=1}^N 2^{-i/p}p_i\Big) \Big|\sum_{i=1}^n
u_{n,i}p_i\Big|^2 \Big(\sum_{i=1}^N 2^{-i/p}p_i\Big) \leq
a^2.$$ 
 We may apply Lemma \ref{lem:counter} to get $\lambda_i\geq 0$ with $\sum_{i=1}^N \lambda_i=1$ and $a=(2 C)^{-1} \sum_{i=1}^N 2^{-i/p}\lambda_i^{-1/2} p_i$ for some contraction $C$.
 
Note that at least one of the $\lambda_i$ is smaller than $\frac 1 N$.  It follows that $\|a\|_{p,\infty}\gtrsim N^{1/2} 2^{- i/p} \|p_i\|_{p,\infty}=N^{1/2}$ whereas $\|X\|_{\Lambda_{p,q}^c}\lesssim N^{1/q}$. 

We also have the estimate $ 2\|a\|_p\geq \big(\sum_{i=1}^N \lambda_i^{-p/2}\big)^{1/p}$. Thanks to the H\"older inequality with $q=1+\frac p 2$:
$$ N=\sum_{i=1}^N \frac{\lambda_i^{p/(2q)}}{ \lambda_i^{p/(2q)}}\leq  \Big(\sum_{i=1}^N  \lambda_i^{-p/2}\Big)^{1/q} \Big(\sum_{i=1}^N {\lambda_i}\Big)^{1/q'}.$$
Thus, similarly $2\|a\|_p\geq N^{1/p+1/2}$ whereas $\|X\|_{\Lambda_{p}^c}\lesssim_p N^{1/p}$.
\end{proof}
\begin{rk}\label{amel}
  { \rm The same proposition applies to row spaces by taking adjoints. Actually, it is easy to see that the formal
    identity map on $L_p(\ell_\infty^c)\cap L_\infty(\ell_\infty^c)$ is not
    bounded from
    $\Lambda_{p,q}$ to
    $L_{p,\infty}(\ell_\infty^c)+L_{p,\infty}(\ell_\infty^r)$ using the same counter-example. Indeed, let $s\in \mathbb
B(\ell_2)$, be the shift operator: if $(e_i)_{i\geq 0}$ is the canonical basis of $\ell_2$ then $se_i = e_{i+1}$ for every $i\geq 0$. We can assume that the set
$\{(u_{n,i});n\geq 0 \}$ is also stable by left multiplication by
$s$. Now, assume that $X=(x_n)_n$ in the proof above decomposes as $Y+Z$ with $Y=(w_n)a$,
  $Z=b(v_n)$
with $(v_n), (w_n)$ bounded and $a,b\in L_{p,\infty}$. For 
any $k\geq 0$, by assumption on $(u_{n,i})$, $s^kX$ is a subsequence of $X$ so there are indices $n_k$ such
that $S^kx_n=w_{n_k}a+ bv_{n_k}$. Thus $x_n=(S^*)^kw_{n_k}a+ (S^*)^kbv_{n_k}$ for all $n$.
Since $b$ is in the Schatten $(p,\infty)$-class,  $S^{*k}b$ goes to 0 in $(p,\infty)$-norm and we can assume by weak-$*$ compacity that $x_n=\tilde w_n a$ for some bounded  
 sequence $(\tilde w_n)$. This implies that, $\|X\|_{L_p(\ell_\infty^c)+L_p(\ell_\infty^r)}=\|X\|_{L_p(\ell_\infty^c)}$.
 }
\end{rk}

\begin{rk}{\rm When the second index $q<2$ and $2< p<\infty$,
    we have an inclusion
    $$\Lambda_{p,q}(\ell_\infty^c)\subset L_{p,\infty}(\ell_\infty^c).$$
One can see it using duality. The space
 $L_{p,q}(\ell_\infty^c)$ is the anti-dual of 
 $L_{p',q'}(\ell_1^c)$ when $q\geq 2$ with $\frac 1 p+ \frac 1{p'}= 1$ and $\frac 1 q+ \frac 1{q'}= 1$. It consists of sequences $(x_n)$ such that there are 
column contractions $(u_n)$ and $(v_n)$ in $\N$ and $\alpha\in L_2(\N)$
and $\beta\in L_{p^\#,q^\#}(\N)$ where $1/2+1/q^\#=1/q'$ and $1/2+1/p^\#=1/p'$ with 
$x_n=\alpha u_n^*v_n\beta$. The norm is obtained by taking the infimum of 
$\|\alpha\|_2.\|\beta\|_{p^\#,q^\#}$ over all possible decompositions.  

Using 
$L_{p^\#,q^\#}=(L_{2^\#},L_{2})_{1- \frac 2p,q^\#}$, it follows that the
inclusion $L_{p',q'}(\ell_1^c)\subset (
L_{2}(\ell_1^c),L_{1}(\ell_1^c))_{1-\frac 2 p,q^\#}$ is bounded. This yields a
bounded inclusion 
$$(L_{2}(\ell_\infty^c),L_\infty(\ell_\infty^c))_{1-\frac 2p,\frac {2q}{q+2}}\subset L_{p,q}(\ell_\infty^c).$$
In particular,  $\Lambda_{\theta,q}(\ell_\infty^c)\subset \Lambda_{\theta,2}(\ell_\infty^c)\subset L_{p,\infty}(\ell_\infty^c)$.
}
\end{rk}

The version of the Marcinkiewicz theorem for $\Lambda$-spaces can now be written as 

\begin{cor}\label{interpo}
  Assume that $S$ is of weak type $(p_0,p_0)$ and $(p_1,p_1)$ with $0<p_0<p_1\leq \infty$. Let $p_0<p<p_1$, then for any $X\in L_p(\M)$, there exist $Z\in \Lambda_p(\ell_\infty^c)$, $Y\in \Lambda_p(\ell_\infty^r)$ so that $S(x)=Z+Y$ and 
$$
\|Z\|_{\Lambda_p^c}+\|Y\|_{\Lambda_p^r} \leq C_{p} \|X\|_p.
$$
\end{cor}
Actually, $Z$ and $Y$ only depend on $X$, not on $p$.

\begin{rk}{\rm The previous arguments can be used to actually prove
    that one can not replace $\Lambda_p^c+\Lambda_p^r$ by
    $L_p(\ell_\infty^c)+L_p(\ell_\infty^r)$ in this corollary 
    when $0<p_1\leq\infty$ in general. Indeed, using the notation from Proposition
    \ref{Ll}, set $U_i=(u_{n,i})_{i\geq 0}$ and define a map
    $S: \ell_0^{\omega^2}\to B(\ell_2)^{\bN}$ by
    $S((x_i))=\sum_{i\geq0}^N x_iU_ip_i$. By Lemma \ref{decomp}, it is
    bounded from $\ell_q^{\omega^2}$ to $\Lambda_q^c$ for all
    $q<\infty$  and thus of weak type $(q,q)$ with constant $C_q$ independent of $N$. For $q=\infty$, its norm
    is controlled by $\sqrt N$.
    Taking $x=(2^{-2i/p})$ as in \ref{Ll}, we have $\|x\|_p=N^{1/p}$ but Remark
    \ref{amel} yields
    $$\|S(x)\|_{L_p(\ell_\infty^c)+L_p(\ell_\infty^r)}=\|S(x)\|_{L_p(\ell_\infty^c)}\gtrsim N^{1/2+1/p}.$$
 This  covers the case $p_1<\infty$ choosing $p\in (\max\{p_0,2\},p_1)$. For $p_1=\infty$, one just need to note that homogeneity would imply that for such $p$ and $\theta\in (0,1)$ such that $\frac 1 p =\frac {1-\theta}{p_0}$: $$\sqrt N\lesssim \|S\|_{ \ell_p^{\omega^2}\to L_p(\ell_\infty^c)+L_p(\ell_\infty^r) }\lesssim  C_{p_0}^{1-\theta}{\sqrt N }^{\theta}.$$
}\end{rk}

Let us conclude with two variations of the previous arguments involving $L_p(\ell_\infty^c)+L_p(\ell_\infty^r)$. The first is about projections.

\begin{lemma}\label{r+c}
  Assume that $S$ is of weak type $(p_0,p_0)$ and $(p_1,p_1)$ with
  constants $1$ with $p_0, p_1\geq 1$.  Then for
  any finite projection $r\in \M$, there exist disjoint projections $q_k$, $k\in \bZ$, with $\tau(q_k)\leq 2^k\tau(r)$ sequences $u_n,v_n\in \N$ and a
  constant $C=C_{p_0,p_1}>0$
  $$\forall \;n\geq 0,\;   S_n(r)=z u_n+ v_n z , \quad \textrm{ with }
  z=\sum_{k}  c_k q_k, \quad \textrm{and }  \|u_n\|,\|v_n\|\leq C,$$
  where $c_k=2^{k/p_1} (|k|+1)$ for $k\leq0$ and $c_k=2^{-k/p_0} k$  for $k \geq 1$.

  In particular, for $p\in (p_0,p_1)$ 
   $$
   \norm{S(r)}_{L_p(\ell_\infty^c) + L_p(\ell_\infty^r)} \leq C_{p,p_0,p_1} \norm{r}_p.
   $$
\end{lemma}
\begin{proof}
With the notation from the proof of Lemma \ref{basic}, 
$$S_n(r)= \sum_{k\in \bZ} \big(q_kS_n(r) e_{k+1}+ e_kS_n(r) q_k\big).$$
For $k\in \bZ$, set $b_k = \min (2^{k/p_0},2^{k/p_1})$ and put
$$
  u_n=\sum_{k\in \bZ} \frac{b_k} {\md{k}+1} q_kS_n(r) e_{k+1} \quad v_n=\sum_{k\in \bZ} \frac{b_k} {\md{k}+1} e_kS_n(r) q_{k} \quad
 z=\sum_{k\in \bZ} b_k (\md{k}+1) q_k.
$$
Then one easily checks $S_n(r)=z u_n+ v_n z$.
$\|u_n\|, \|v_n\|\leq C_{p_0,p_1}$ and the estimate for $z$ as well as
$S_n(r)=zu_n+v_nz$. The last estimate is also clear.
\end{proof}

\begin{rk}{\rm
For $a\in L_{2p}^+$, $p\geq 1$ and any contraction $u\in \M$, one can
always find a contraction $v\in \M$ such that $aua=\frac 12 (a^2v+va^2)$;
this follows from instance from the Cauchy formula for the holomorphic function 
$F(z)=a^{2(1-z)}ua^{2z}$ with values in $L_p$. Thus, when $S$ is positive, Lemma \ref{basic} implies this one.}
 \end{rk}

For general elements, we can get

\begin{prop}\label{r+c2}
  Assume that $S$ is of weak type $(p_0,p_0)$ and $(p_1,p_1)$ with constants 1.
  Then for any $x\in L_{p_0}(\M)\cap L_{p_1}(\M)$, there exist $z\in
  \cap_{p_0< p<p_1} L_p(\N)^+$ and sequences $u_n,v_n\in \N$ such that
  $$\forall \;n\geq 0,\;   S_n(x)=z u_n+ v_n z ,\;  \|u_n\|,\|v_n\|\leq 1 \; \textrm{and }\;  \forall\, p_0<p<p_1,\;  \|z\|_p\leq C_{p} \|x\|_{p_0}^{1-\theta}\|x\|_{p_1}^{\theta},$$
  where $\frac 1 p= \frac {1-\theta}{p_0}+\frac \theta {p_1}.$
\end{prop}

\begin{proof}
  Decomposing $x$ as a combination of four positive elements and using the fact that $L_p(\N;\ell_\infty^c)$ and $L_p(\N;\ell_\infty^r)$ are quasi-Banach spaces, we can assume that $x\in L_p^+$.

Using homogeneity, we may assume that $\|x\|_p=1$. 
We use Lemma \ref{decdya} to write $x=\sum_{m\in \bZ} 2^{-m}r_m$ for some finite projections $r_m$.

We can apply Lemma \ref{r+c} to each $r_m$ to get
$$S_n(x)= \sum_{m\in \bZ} 2^{-m} (z_m u_{m,n} + v_{m,n}z_m),$$
with $\|u_{m,n}\|,\|v_{m,n}\|\leq 1$ and $\mu(z_m)\preceq \sum_{k\in \bZ}  2^{-k} c_kD_{2^k}(\mu(r_m))$. 

Using Lemma \ref{fact}, we can find elements $u_n$ and $v_n$ in $\N$ such that
$$ \sum_{m\in \bZ} 2^{-m} z_m u_{m,n} =  \Big( \sum_{m\in \bZ} 2^{-2m} m^2 z_m^2\Big)^{1/2}u_n,\quad \sum_{m\in \bZ} 2^{-m}v_{m,n} z_m =  v_n\Big( \sum_{m\in \bZ} 2^{-2m} m^2 z_m^2\Big)^{1/2},$$
with  $\|u_n\|, \|v_n\|$ bounded by an absolute constant $C$.

Set $z=\Big( \sum_{m\in \bZ } 2^{-2m} m^2 z_m^2\Big)^{1/2}$. Then $z \otimes e_{0,0}$ is 
the modulus of $\sum_{m\in \bZ} 2^{-m} m z_m\otimes e_{m,0}$ in 
$L_p(\N \overline{\otimes} \B(\ell_2))$. Thus
$$ \mu(z) \preceq \sum_{\in \bZ}  2^{-m}m \sum_{k\in \bZ} c_k D_{2^k}(\mu(r_m))=
\sum_{k\in \bZ} c_k D_{2^k}\big(\sum_{m\in \bZ} 2^{-m} m \mu(r_m)\big).$$

By Lemma \ref{decdya},
$ \sum_{m>0} 2^{-m} m \mu(r_m)\leq C \mu\big(x (|\ln(x)|+1)\big)$ for
some constant $C>0$.

We get $$\|z\|_p\lesssim \sum_{k\in \bZ} c_k2^{k/p} \|x (|\ln(x)|+1)\|_p.$$
where we used that $D_{2^k}$ has norm less than $2^{k/p}$ on $L_p$.

We have the inequalities $t^p (|\ln(t)|+1)^p\lesssim_{p_,p_0}t^{p_0}$
for $0<t<1$ and $t^p (|\ln(t)|+1)^p\lesssim_{p_,p_1}t^{p_1}$ for
$t\geq 1$. They yield that $\|x (|\ln(x)|+1)\|_p^p\lesssim \|x\|_{p_0}^{p_0}+\|x\|_{q_1}^{q_1}$.

Hence we have that
$\|z\|_p\lesssim \|x\|_{p_0}^{p_0}+\|x\|_{q_1}^{q_1}$.  Using

The conclusion follows using homogeneity as usual.
\end{proof}
Unfortunately, trying to improve the above estimate using reiteration as in \cite{JX2} does not provide anything more than Corollary \ref{interpo}.

\bigskip

Let us conclude by highlighting the link between $\Lambda$-spaces and bilateral almost uniform convergence. Recall that a sequence of operators $(x_n)$ converges {\it bilaterally almost uniformly} (b.a.u.) to $x$ if and only if for any $\varepsilon>0$, there exists a projection $e\in \M$ such that 
$$
\tau (1-e) \leq \varepsilon \quad \textrm{and} \quad
\norm{e(x_n - x)e} \to 0.
$$
This notion was introduced to serve as one of the non commutative analogs of pointwise convergence (see \cite{Jaj85} for more details and the early developments related to this form of convergence). Let $(S_n)_{n\geq 0}$ be a sequence of operators defined on $L_1(\M) + \M$. A standard approach to prove the convergence of $(S_n)_{n\geq 0}$ to some limit operator $S_\infty$ on $L_p(\M)$ is the following : first, prove convergence on a dense subset of $L_p(\M)$ (usually taken to be $L_1(\M) \cap \M$); second, prove a maximal inequality and use a form of Banach principle to extend the convergence on all of $L_p(\M)$. 

\begin{prop}
 Let $p\in [1,\infty)$ and $S = (S_n)$ be a sequence of maps from $L_p(\M)$ to $L_p(\N)$ such that $S$ is bounded from $L_p(\M)$ to $\Lambda_p(\N;\ell_\infty)$. Assume that there exists a bounded linear map $S_\infty$ from $L_p(\M)$ to $L_p(\N)$ and a dense subset $E \subset L_p(\M)$ such that for any $x\in E$, $(S_n(x))$  converges b.a.u. to $S_\infty(x)$. Then $(S_n(x))$ converges b.a.u. to $S_\infty(x)$ for any $x\in L_p(\M)$.
\end{prop}

As the topic of b.a.u. convergence is not central to this paper, the proof is left to the reader. The intended application for such a result would be to the convergence of sequences of operators that are not positive, verify weak type maximal inequalities and for which the strong type cannot be proved. But no natural example of such a behaviour has been discovered so far as every proof of maximal inequality relies in one way or another on a reduction to positive operators. A good candidate would be the matrix valued Carleson operator. In the meantime, we believe that the proposition (or more refined versions that would involve symmetric spaces) complements results of b.a.u. convergence in tracial von Neumann algebras in the spirit of \cite{CL21}. 

\bigskip

\textbf{Acknowledgment.}  The authors are supported by
ANR-19-CE40-0002. The authors thanks the Lorentz Center in which discussions 
about this subject were initiated. 

\bibliographystyle{plain}

\bibliography{bibli}

\end{document}